\documentclass[a4paper, 11pt]{article}
\usepackage{amsmath, amssymb, amsthm, latexsym, enumerate}
\usepackage[alphabetic]{amsrefs}
\usepackage{xspace}
\usepackage{a4wide}
\usepackage[usenames,dvipsnames,svgnames]{xcolor}
\usepackage{tikz}
\usepackage{hyperref}
\usepackage{authblk}
\usepackage{chngcntr}
\usepackage{enumitem}

\newtheorem*{ctheorem}{Theorem \hyperref[bthm:key]{\ref{bthm:key}.\ref{key:quasiisom}'}}

\newtheorem{proposition}{Proposition}[section]
\newtheorem{theorem}[proposition]{Theorem}
\newtheorem{corollary}[proposition]{Corollary}

\newtheorem*{theorem*}{Theorem}

\theoremstyle{definition}
\newtheorem{definition}[proposition]{Definition}
\newtheorem{question}[proposition]{Question}
\newtheorem{remark}[proposition]{Remark}
\newtheorem{example}[proposition]{Example}

\newcommand{\set}[1]{\left\{#1\right\}}
\newcommand{\setcon}[2]{\left\{#1\ \left|\ #2\right.\right\}}

\newcommand{\into}{\hookrightarrow}
\newcommand{\C}{\mathbb{C}}
\newcommand{\R}{\mathbb{R}}
\newcommand{\Z}{\mathbb{Z}}
\newcommand{\N}{\mathbb{N}}

\newcommand{\fgen}[1]{\left\langle #1 \right\rangle}

\newcommand{\fpres}[2]{\left\langle #1 \left| #2 \right.\right\rangle}

\newcommand{\cdim}[1]{\textrm{cdim}(#1)}

\newcommand{\Mcd}{\textrm{cdim}_{\partial_M}}
\newcommand{\stabdim}{\textrm{asdim}_s}

\newcommand{\Mod}{\mathrm{MCG}(\Sigma)}
\newcommand{\teich}{\mathcal{T}(\Sigma)}
\newcommand{\PSL}{\mathrm{PSL}}
\newcommand{\Out}{\mathrm{Out}(\mathbb F)}

\title{A survey on Morse boundaries \& stability}

\author{Matthew Cordes}

\date{\today}

\begin{document}

\maketitle

\section{generalizing hyperbolicity}


Recently there have been efforts to generalize tools used in the setting of Gromov hyperbolic spaces to larger classes of spaces. We survey here two aspects of these efforts: Morse boundaries and stable subspaces. Both the Morse boundary and stable subspaces are systematic approaches to collect and study the hyperbolic aspects of finitely generated groups. We will see in this survey a nice relationship between the two approaches. 

Throughout the survey we will expect the reader to be familiar with the basics of geometric group theory. For details on this material see \cite{bridson:1999fj} and for additional information especially on boundaries of not-necessarily-geodesic hyperbolic spaces see \cite{Buyalo:2007ab}.

We begin with an important definition that has its roots in a classical paper of Morse \cite{morse:1924aa}:

\begin{definition} A geodesic $\gamma$ in a metric space is called \emph{$N$-Morse} if there exists a function 
$N=N(\lambda, \epsilon)$ such that for any $(\lambda, \epsilon)$-quasi-geodesic $\sigma$ with endpoints on $\gamma$, we have $\sigma \subset \mathcal{N}_N(\gamma)$, the $N$-neighborhood of $\gamma$. We call the function $N\colon\R_{\geq1} \times \R_{\geq0}\to\R_{\geq0}$ a \emph{Morse gauge}. A geodesic is Morse if it is $N$-Morse for some $N$.
\end{definition} 

In a geodesic $\delta$-hyperbolic space, the Morse lemma says there is a Morse gauge $N$ which depends only on $\delta$ such that every ray is $N$-Morse, i.e., the behavior of quasi-geodesics (on a large scale) is similar to that of geodesics. This property fails in a space like $\R^2$. On the other hand, if every geodesic in some geodesic space is $N$-Morse, then the space is $\delta$-hyperbolic, where $\delta$ depends on $N$. 

There are a few other competing definitions of geodesics which admit  ``hyperbolic like" properties: geodesics which satisfy a contracting property,  geodesics with superlinear divergence, and having cut points in the asymptotic cone. (See Section \ref{sec: contracting} for the definition of the contracting property and \cite{charney:2015aa, Arzhantseva:2016aa, drutu:2005aa} for the others.) All of these notions have been used extensively to analyze many groups and spaces: right-angled Artin groups \cite{behrstock:2012aa, koberda:2014aa, cordes:2016aa}, Teichm\"uller space \cite{minsky:1996aa, Behrstock:2006aa, Sultan:2014aa, Brock:2006aa, Brock:2011aa}, the mapping class group and curve complex {\cite{Masur:1999aa, Masur:2000aa, durham:2015aa, bestvina:2015ab}},  $\mathrm{CAT}(0)$ spaces \cite{Behrstock:2014aa, Sultan:2014aa}, $\mathrm{Out}(F_n)$ \cite{Algom-Kfir:2011aa}, relatively hyperbolic groups and spaces \cite{drutu:2005aa, Osin:2006aa},  acylindrically hyperbolic groups \cite{dahmani:2011aa, sisto:2016aa, bestvina:2015ab}, and small cancellation groups \cite{Arzhantseva:2016ab}  among others. Furthermore these properties find applications in rigidity theorems such as Mostow Rigidity in rank-1 \cite{Paulin:1996aa} and the Rank Rigidity Conjecture for $\mathrm{CAT}(0)$ spaces \cite{Ballmann:1995aa, Bestvina:2009aa, Caprace:2010aa}.

Geodesics of these types all  have a relationship to the Morse property.  In \cite{charney:2015aa} Charney and Sultan show that Morse, superlinear `lower divergence', and `strongly' contracting are equivalent notions in $\mathrm{CAT}(0)$ spaces. In \cite{Arzhantseva:2016aa} the authors characterize Morse quasi-geodesics in arbitrary geodesic metric spaces and show they are `sublinearly' contracting and have `completely superlinear' divergence. In \cite{drutu:2010aa} the authors show that a quasi-geodesic is Morse if and only if it is a cut point in every asymptotic cone. This plays a key role in showing that metric relative hyperbolicity is preserved by quasi-isometry \cite{drutu:2005aa}. 

While groups may have no Morse geodesics, Sisto in \cite{sisto:2016aa} shows that every acylindrically hyperbolic group has a bi-infinite Morse geodesic. This class of groups, recently unified by Osin in \cite{osin:2016aa}, encompasses many groups of significant interest in geometric group theory: hyperbolic and relatively hyperbolic groups, non-directly decomposable right-angled Artin groups, mapping class groups, and $\mathrm{Out}(F_n)$. There are groups outside this class which contain Morse geodesics. Examples of these groups appear in \cite{olcprimeshanskii:2009aa}; these groups have been shown to have infinitely many Morse geodesics but are not acylindrically hyperbolic. At the other extreme, any groups which admit a law do not contain any Morse geodesics \cite{drutu:2005aa}. It is easy to see that $\R^n$ for $n \geq 2$ has no Morse geodesics: take any geodesic ray $\gamma$ and consider a right angled isosceles triangle whose hypotenuse is on $\gamma$. It is easy to check that the other two edges form a $(\sqrt{2}, 0)$-quasi-geodesic with endpoints on $\gamma$. This is true  no matter how large the triangle and thus these uniform quasi-geodesics can be arbitrarily far from $\gamma$ and therefore violate the Morse condition. 
\medskip

In this survey we will overview the notion of Morse boundaries with increasing generality in Sections \ref{sec: contracting & morse} and \ref{sec:Metric Morse}. In Section \ref{sec:Metric Morse} we will also introduce stability, the notion of stable equivalence, two `stable equivalence' invariants and end with some calculations of these invariants. In Section \ref{sec:stable subgroups} we will spend some time discussing stable subgroups and introduce a boundary characterization of convex cocompactness which is equivalent to stability. Finally, in Section \ref{sec: Metrizing} we will discuss a topology on the contracting boundary which is second countable and thus metrizable. 

\section{contracting and Morse boundaries} \label{sec: contracting & morse}

Boundaries have been an extremely fruitful tool in the study of hyperbolic groups. The visual boundary, as a set, is made up of equivalence classes of geodesic rays, where one ray is equivalent to the other if they fellow travel. Roughly, one topologizes the boundary by declaring open neighborhoods of a ray $\gamma$ to be the rays that stay close to $\gamma$ for a long time. Gromov showed that  a quasi-isometry between hyperbolic metric spaces $X$ and $Y$ induces a homeomorphism on the visual boundaries. In the setting of a finitely generated group $G$ acting geometrically on $X$ and $Y$, then the quasi-isometry from $X$ to $Y$ induced by these actions extends $G$-equivariantly to a homeomorphism of their boundaries. This means, in particular, that the boundary of a hyperbolic group (as a topological space) is independent of the choice of (finite) generating set.

The boundary of a hyperbolic group is a powerful tool to study the structure of the group. For instance, Bowditch and Swarup relate topological properties of the boundary to the JSJ decomposition of a hyperbolic group \cite{Bowditch:1998aa, Swarup:1996aa}. Bestvina and Mess in \cite{Bestvina:1991aa} relate the virtual cohomological dimension of a hyperbolic group $G$ to the dimension of the boundary of $G$. Boundaries have also been instrumental in the proofs of quasi-isometric rigidity theorems, particularly Mostow Rigidity in rank 1 \cite{Paulin:1996aa}.

There is also a robust boundary theory of relatively hyperbolic groups using the Bowditch boundary introduced in a 1998 preprint of Bowditch that was recently published \cite{Bowditch:2012aa}. Bowditch used this boundary to analyze JSJ decompositions \cite{Bowditch:1998aa, Bowditch:2001aa}. Groff showed that if a group $G$ is hyperbolic relative to a collection $\mathcal{A}$ of subgroups which are not properly relatively hyperbolic, then the Bowditch boundary is a quasi-isometry invariant \cite{Groff:2013aa}.

The topology on the visual boundary for a $\mathrm{CAT}(0)$ space can be defined in a manner similar to the visual boundary of a hyperbolic space. Hruska and Kleiner show that in the case of a $\mathrm{CAT}(0)$ space with isolated flats, this boundary is a quasi-isometry invariant \cite{Hruska:2005aa}. Unfortunately, Croke and Kleiner produced an example of a right-angled Artin group which acts geometrically on two quasi-isometric $\mathrm{CAT}(0)$ spaces which have non-homeomorphic visual boundaries \cite{Croke:2000aa}, hence the boundary of a $\mathrm{CAT}(0)$ group is not well-defined. More surprising is that Wilson shows that this group admits an uncountable collection of distinct boundaries \cite{Wilson:2005aa}.  Charney and Sultan in \cite{charney:2013qy} showed that if one restricts attention to rays with hyperbolic-like behavior, \emph{contracting rays}, then one can construct a quasi-isometry invariant boundary for any complete $\mathrm{CAT}(0)$ space. They call this boundary the contracting boundary. 

\subsection{contracting boundaries} \label{sec: contracting}

We begin with a description of the contracting boundary of a $\mathrm{CAT}(0)$ space introduced in \cite{charney:2015aa}. We will start with some explication on contracting geodesics.

\begin{definition}[contracting geodesics] \label{def: str contracting} Given a fixed constant $D$, a geodesic $\gamma$ is said to be \emph{$D$-contracting} if for all $x, y \in X$, $$d_X(x,y) < d_X(x, \pi_\gamma(x)) \implies \mathrm{diam}(\pi_\gamma(x), \pi_\gamma(y))<D$$ where $\pi_\gamma$ is the closest point projection to $\gamma$. We say that $\gamma$ is \emph{contracting} if it is $D$-contracting for some $D$. An equivalent definition is that any metric ball $B$ not intersecting $\gamma$ projects to a segment of length $<2D$ on $\gamma$.
\end{definition}

\begin{remark}This definition is sometimes referred to as strongly contracting. To see a more general definition of contracting and its relationship with the Morse property and divergence see Section \ref{sec: Metrizing} and for more details \cite{Arzhantseva:2016aa}. This more general definition is used by Cashen to answer a question of Charney--Sultan (see Theorem \ref{thm:cashen}) and by Cashen--Mackay to construct a topology on the contracting boundary that is metrizable (see Section \ref{sec: Metrizing}). \end{remark}

 It is easy to see that any geodesic in a flat will not be contracting because the projection of any ball onto a geodesic will just be the diameter of the ball. In a hyperbolic $\mathrm{CAT}(0)$ space all geodesics are uniformly contracting. Contracting geodesics appear in many groups and spaces, for instance: psuedo-Anosov axes in Teichm\"uller space \cite{minsky:1996aa}, iwip axes in the Outer Space of outer automorphisms of a free group \cite{Algom-Kfir:2011aa}, and axes of rank 1 isometries of $\mathrm{CAT}(0)$ spaces \cite{Ballmann:1995aa, Bestvina:2009aa}. 
 
\begin{remark}[Morse and (strongly) contracting are not equivalent] \label{rem: contracting neq morse} If $X$ is proper geodesic space, then contracting implies Morse \cite{Algom-Kfir:2011aa}. In this generality the converse is not true. Consider a ray $Y$ and a set of intervals $\{I_i\}_{i \in \N}$ of length $i$ on $Y$. To the endpoints of these intervals attach an edge of length $i^2$. Let $X$ be the resulting space. The ray $Y \subset X$ is not contracting because balls can have unbounded projection onto $Y$. To see that $Y$ is Morse, we need to check to check that there exits an constant $N=N(\lambda, \epsilon)$ so that for any $(\lambda, \epsilon)$-quasi-geodesic $\sigma$ with endpoints on $Y$, $\sigma \subset \mathcal{N}_N(Y)$. By \cite[III.H Lemma 1.11]{bridson:1999fj} we can replace $\sigma$ with a `tame' quasi-geodesic $\sigma'$ with three nice properties: $\sigma$ and $\sigma'$ have the same endpoints, $\sigma$ and $\sigma'$ have Hausdorff distance less than $(\lambda+\epsilon)$, and the length between any two points $x,y$ on of $\sigma'$ is bounded by a linear function of the $d(x,y)$ (with constants depending only on $\lambda, \epsilon$). Since the length of the attached intervals is $i^2$ while the distance between their endpoints is $i$, in order for $\sigma'$ to cross such a segment, we must have $i^2$ is less than a linear function of $i$. Thus $\sigma'$ must only cross finitely many of the attached segments. Thus for for any $(\lambda,\epsilon)$-quasi-geodesic with endpoints on $Y$ we can choose that $N(\lambda, \epsilon)=j^2+ \lambda +\epsilon$ where $j^2$ is the length of the longest attached interval $\sigma'$ travels over.

In the setting of $\mathrm{CAT}(0)$ spaces Sultan shows that Morse and contracting are equivalent notions \cite{Sultan:2014aa} and in \cite{charney:2015aa} Charney and Sultan reprove this fact with explicit control on the constants.  Charney and Sultan also characterize contracting geodesics in $\mathrm{CAT}(0)$ cube complexes using a combinatorial criterion which gives an effective tool for analyzing the boundary. We will use this characterization to help us compute Examples \ref{example:z fp zsqd} and \ref{example: croke--kleiner}. \end{remark}

Let $X$ be a $\mathrm{CAT}(0)$ space. We define the visual boundary, $\partial X$, to be the set of equivalence classes of geodesic rays up to asymptotic equivalence and denote the equivalence class of a ray by $\alpha(\infty)$. It is an elementary fact that, for $X$ a complete $\mathrm{CAT}(0)$ space and $e\in X$ a fixed basepoint, every equivalence class can be represented by a unique geodesic ray emanating from $e$. One natural topology on $\partial X$ is the cone topology. We define the topology of the boundary with a system of neighborhood bases. A neighborhood basis for $\alpha(\infty)$ is given by open sets of the form:
 \[U(\alpha, r, \epsilon) = \{\beta(\infty) \in \partial X \mid \beta \text{ is a geodesic ray based at } e \text{ and } \forall t<r, d(\beta(t), \alpha(t))< \epsilon \}.\] That is, two geodesic rays in the cone topology are close if they fellow travel (at distance less than $\epsilon$) for a long time (at least time $r$). This topology is independent of choice of basepoint. When we refer to the visual boundary we will always assume that means with the cone topology unless otherwise stated.

Let $X$ be a complete $\mathrm{CAT}(0)$ space with basepoint $p \in X$. We define the \emph{contracting boundary} of a $\mathrm{CAT}(0)$ space $X$ to be the subset of the visual boundary consisting of all contracting geodesics: \[\partial_c X_e = \{ \alpha(\infty) \in \partial X \mid \alpha \text{ is contracting with basepoint } e \}.\]

In order to topologize the contracting boundary we consider a collection of increasing subsets of the boundary, \[\partial_c^n X_e = \{ \gamma(\infty) \in \partial X \mid \gamma(0)=e, \gamma \text{ is a $n$-contracting ray}\},\] one for each $n \in \N$. We topologize each $\partial_c^n X_e$ with the subspace topology from the visual boundary of $X$. We note that there is an obvious continous inclusion map $i \colon \partial_c^m X_e \into \partial_c^n X_e$ for all $m<n$. We can topologize the whole boundary by taking the direct limit over these subspaces. Thus $\partial_c X_e=\varinjlim \partial_c^n X_e$ with the direct limit topology. Recall that this means a set $U$ is open (resp. closed) in $\partial_c X_e$ if and only if $U \cap \partial_c^n X_e$ is open (resp. closed) for all $n \in \N$.

One nice property of this boundary is that if we fix a contracting constant $n$, then $\partial_c^n X_e$ is compact and behaves very much like the boundary of a hyperbolic space, an analogy we will build on in Section \ref{sec:Metric Morse}.

Another nice property is that the direct limit topology on $\partial_c X_e$ does not depend on the choice of basepoint so we will freely denote the contracting boundary as $\partial_c X$ without mention of basepoint. 

\subsubsection{examples}
We will focus on examples for the contracting boundary for a couple reasons. First, other more general incarnations of the Morse boundary, which we will see later, are homeomorphic to the contracting boundary in the case of a $\mathrm{CAT}(0)$ space. Second, there is a combinatorial characterization of contracting geodesics in $\mathrm{CAT}(0)$ cube complexes (see \cite[Theorem 4.2]{charney:2015aa}) that makes the computations more transparent.

\begin{example}[$\R^n$, hyperbolic spaces]Revisiting the conversation of which spaces have contracting geodesics, when $X= \R^n$ for $n \geq 2$ or more generally the product of two unbounded $\mathrm{CAT}(0)$ space will have no contracting geodesics because any geodesic will be contained in a flat. This means that $\partial_c X$ will be empty. On the other hand, if $X$ is $\mathrm{CAT}(0)$ and hyperbolic, then $\partial_c X$ will be the Gromov boundary because every ray will uniformly contracting. For example, $\partial \mathbb{H}^n$ will be homeomorphic to $\mathbb{S}^{n-1}$. \end{example}

\begin{example}[hyperbolic with a dash of flat] Now consider the space $X$ formed by gluing a Euclidean half-plane to a bi-infinite geodesic $\gamma$ in the hyperbolic plane $\mathbb{H}^2$. Picking $\gamma(0)$ to be our basepoint, we see that any geodesic in the half-plane (including $\gamma$) will not be Morse, but all the other rays in $\mathbb{H}^2$ will still be contracting. So the boundary will be the union of two open intervals. We can already see here how the contracting boundary differs from the Gromov boundary, because it is not necessarily compact. In fact, Murray shows that the contracting boundary of a $\mathrm{CAT}(0)$ space is compact if any only if $X$ is hyperbolic \cite{murray:2015aa}. \end{example}

\begin{figure} \label{fig:ZfreeZsqrd}
  \centering
  \def\svgwidth{2in}
  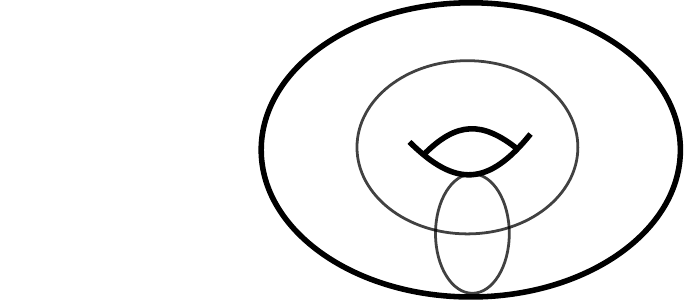
\caption{presentation complex of $\Z \ast \Z^2$}
\end{figure}

\begin{example}[$\Z \ast \Z^2$] \label{example:z fp zsqd}For a slightly more complicated example consider the space $Y$ formed by the wedge product of a loop and a torus and its fundamental group $\pi_1(Y)=\Z \ast \Z^2=\langle a \rangle \ast \langle b, c \rangle$. See Figure \ref{fig:ZfreeZsqrd}. The universal cover of $Y$, $\tilde{Y}$, is a $\mathrm{CAT}(0)$ cube complex \cite{Charney:1995aa}. Consider a geodesic ray $\gamma$ in $X$. It is not hard to see that the contracting constant of $\gamma$ only depends on how long subsegments of $\gamma$ spend in the cosets of the $\Z^2$ subgroup (because of the tree-graded structure of $\tilde{Y}$), and it is only contracting if there is a uniform bound on the length of these subsegments. In fact, the contracting constant is half that uniform bound; thus each $\partial_c^n \tilde{Y}$ is the set of geodesic rays which have subsegments of lenght at most $2n$ in the $\Z^2$ cosets with the subspace topology. So $\partial_c \tilde{Y}$ is the direct limit over these spaces. Later we will see in Theorem \ref{thm:RAAG} that the subspace of $\tilde{Y}$ formed by the union of all $n$-contracting geodesics with basepoint $e$ is quasi-isometric to a proper simplicial tree and thus has boundary homeomorphic to a Cantor set.\end{example}

\begin{remark}[the contracting boundary is not in general first-countable] \label{rem: not 1st countable} The topology of the contracting boundary is in general quite fine relative to the subspace topology. In fact, in \cite{murray:2015aa} Murray shows that the space $\tilde{Y}$ as defined above is not first-countable (and thus not metrizable). We present Murray's argument here: We choose the basepoint of $\tilde{Y}$ to be some lift of the wedge point in $Y$.  We first note that the geodesic ray $\alpha = aaaa\ldots$ is $0$-contracting. We define a collection of geodesic rays $\beta_i^j=a^{i}b^jaaaa\ldots$. Note that the each of the $\beta_i^j$ are $j$-contracting and not $j'$-contracting for any $j'<j$. If we fix a $j$, it is clear that that the $\{\beta_i^j\}$ converge to $\alpha$ in the contracting boundary. Consider a new sequence $\{\beta_{f(j)}^j\}$ where $f$ is any function $f \colon \N \to \N$. We note that the intersection of $\{\beta_{f(j)}^j\}$ with $\partial_c^n \tilde{Y}$ is finite for each $n \in \N$ and thus closed in the subspace topology and therefore in $\partial_c X$. Ergo, the $\{\beta_{f(j)}^j\}$ {cannot} converge to $\alpha$ in the contracting boundary. 

A general fact for all first-countable spaces is that if you have a countable collection of sequences which all converge to the same point, it is always possible to pick a `diagonal' sequence which also converges to that point. That is, if $\partial_c \tilde{Y}$ were first-countable, since $\{\beta_i^j\}$ converges to $\alpha$ for fixed $j$, then there would be a function $f \colon \N \to \N$ so that $\{\beta_{f(j)}^j\}$  converges to $\alpha$. As we saw above, this cannot happen. Since all metric spaces are first countable, this means that the contracting boundary is not, in general, metrizable! In fact, Murray, in the same paper, shows if $X$ is a $\mathrm{CAT}(0)$ space with a geometric action, then $\partial_cX$ is metrizable if any only if $X$ is hyperbolic.
\end{remark}

\begin{figure} \label{fig:CK}
  \centering
  \def\svgwidth{3in}
  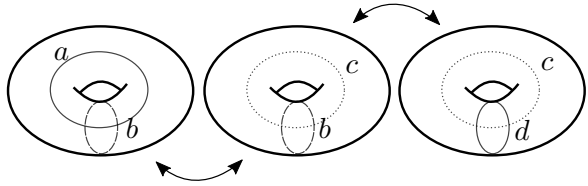
\caption{the Croke--Kleiner example}
\end{figure}

\begin{example}[Croke--Kleiner] \label{example: croke--kleiner} Recall that Croke and Kleiner produced an example of a right-angled Artin group that acts geometrically on two quasi-isometric $\mathrm{CAT}(0)$ spaces which have non-homeomorphic visual boundaries \cite{Croke:2000aa}. Their example is a right-angled Artin group \[A_\Gamma = \fpres{a,b,c,d}{[a,b]=[b,c]=[c,d]=1}.\] The Salvetti complex of this group, $S_\Gamma$, is three tori with the middle torus glued to the other two along orthogonal curves corresponding to the generators $b$ and $c$. See Figure \ref{fig:CK}. As in the example above, it follows from \cite{Charney:1995aa} that the universal cover of $S_\Gamma$, $\tilde{S_\Gamma}$, is a $\mathrm{CAT}(0)$ cube complex. 

Let $B_1$ be the union of the $(a,b)$-torus and the $(b,c)$-torus in $S_\Gamma$ and $\tilde{B_1}$ its inverse image in $\tilde{S_\Gamma}$. Each component of the inverse image decomposes as the direct product of the Cayley graph of the free group on two elements and $\R^2$. Thus the contracting boundary of each component of $\tilde{B_1}$ is empty. The same fact holds for $B_2$, the union of the $(b,c)$-torus and the $(c,d)$-torus in $S_\Gamma$. Croke and Kleiner refer to the components as ``blocks". As in the example above, it follows that in order for a geodesic ray $\gamma$ in $\tilde{S_\Gamma}$ to be contracting, there must be a uniform bound on the length of the subsegments of $\gamma$ which intersect the blocks. The converse of this holds by the combinatorial characterization by Charney and Sultan of contracting geodesics in $\mathrm{CAT}(0)$ cube complexes \cite{charney:2015aa}. Thus, $\gamma$ is contracting if and only if there is a uniform bound on the length of subsegments of $\gamma$ which intersect a single block. We will revisit this in Example \ref{CK redux} after we describe some of the properties of the contracting boundary. \end{example}

\subsubsection{properties of the boundary}

This boundary has many desirable properties. First and foremost it is a quasi-isometry invariant. That is, if $X$ and $Y$ are two complete $\mathrm{CAT}(0)$ spaces and $q\colon X \to Y$ is a quasi-isometry, then $q$ induces a homeomorphism $\partial_q \colon \partial_c X \to \partial_c Y$. Furthermore since each $\partial_c^n X$ is compact, the boundary is $\sigma$-compact. Finally it is a visibility space. In summary:

\begin{theorem}[\cite{charney:2015aa}] Given a complete $\mathrm{CAT}(0)$ space $X$, the contracting boundary, $\partial_c X=\varinjlim \partial_c^n X_e$, equipped with the direct limit topology, is 
\begin{enumerate}
\item independent of choice of basepoint;
\item $\sigma$-compact, i.e., the union of countably many compact sets;
\item a visibility space, i.e., any two points in the contracting boundary can be joined by a bi-infinite contracting geodesic; and
\item a quasi-isometry invariant.
\end{enumerate}
\end{theorem}

\begin{example}[Croke--Kleiner redux] \label{CK redux} 
Croke and Kleiner produce the two spaces on which $A_\Gamma$ acts geometrically by modifying the metric on $\tilde{S_\Gamma}$ in a very simple way: they skew the angles between the $b$ and $c$ curves making  the $(b,c)$-cubes parallelograms. This is enough to change the homeomorphism type of the visual boundary. In fact, Wilson shows that any two distinct angles between the $b$ and $c$ curves produces non-homeomorphic boundaries \cite{Wilson:2005aa}. Qing showed that if you keep the angles $\frac{\pi}{2}$ and change the side lengths of the cubes then the identity map does not induce a homeomorphism on the boundary \cite{Qing:2013aa}. 

In the examples of Croke--Kleiner and Wilson, the parts of the boundary which change the homeomorphism type are the parts that come from the intersection of the blocks. These points do not appear in the contracting boundary. Neither do the points in the Qing example. The parts of the boundary which change are the rays which stay longer and longer time in successive blocks. These examples suggest that the restriction to contracting rays may be optimal if you want a quasi-isometry invariant.
\end{example}

In \cite{charney:2015aa}, Charney and Sultan remarked that a quasi-isometry induces a bijection on the contracting rays and the set of contracting rays could be topologized with just the subspace topology from the visual boundary. They asked if a quasi-isometry would induce a homeomorphism on the boundary with this topology. In \cite{Cashen:2016aa} Cashen answers this question in the negative:

\begin{theorem}[\cite{Cashen:2016aa}] \label{thm:cashen} In general a quasi-isometry will not induce a homeomorphism on the contracting boundary with the subspace topology.\end{theorem}
Cashen's examples, though, are pathological in nature and open questions still remain: If $X, Y$ are $\mathrm{CAT}(0)$ spaces with cocompact isometry groups and if $\phi \colon X \to Y$ is a quasi-isometry, does $\phi$ induce a homeomorphism on the contracting boundary with the subspace topology? Are these two boundaries abstractly homeomorphic?

One can also ask when a homeomorphism between the boundaries of two spaces is induced by a quasi-isometry. In the case of a hyperbolic spaces, Paulin \cite{Paulin:1996aa} gives conditions under which this holds.  Recent work of Charney and Murray \cite{Charney:2017aa} gives some analogous conditions for the contracting boundaries of cocompact CAT(0) spaces.

\subsubsection{dynamics on the contracting boundary} \label{sec:dynamics}
In the realm of hyperbolic groups there is a classical notion of North-South dynamics:
\begin{theorem} If $G$ is a hyperbolic group acting geometrically on a proper geodesic metric space $X$ and if $g$ is an infinite order element, then for all open sets $U$ and $V$ with $g^\infty \in U$ and $g^{-\infty} \in V$ there exists an $n\in\N$ so that $g^nV^c \subset U$. \end{theorem}
It is well known that the classical notion of North-South dynamics of axial isometries on the visual boundary of a $\mathrm{CAT}(0)$ space fails because there are isometries which fix whole flats. Rank-1 isometries, though, do act on the visual boundary with North-South dynamics \cite{Ballmann:1995ab, Hamenstadt:2009aa}. Since rank-1 isometries of  $\mathrm{CAT}(0)$ spaces are contracting, one might hope that North-South dynamics hold on the contracting boundary. This is not the case. Again we will present an example from \cite{murray:2015aa}.

\begin{example}[Example \ref{example:z fp zsqd} revisited]Again, we consider $\Z \ast \Z^2=\langle a \rangle \ast \langle b ,c \rangle$ and the spaces $Y$ and $\tilde{Y}$. Let $\alpha$ be an axis for $a$. Let $\beta_i$ be the geodesic defined by the word $a^{-i}b^iaaaa\ldots$. Again we note that the $\beta_i$ do not converge to $\alpha$ in the contracting boundary because $\{ \beta_i \} \cap \partial_c^D \tilde{Y}$ is finite for every $D$ and thus $\{ \beta_i \}$ is closed in $\partial_c \tilde{Y}$. We note that the set $V=(U(\alpha(-\infty), r, \epsilon) \cap \partial_c \tilde{Y})\smallsetminus \{ \beta_i \}$ is an open set around $\alpha(-\infty)$ but for all $n \in \N$ we have $a^n\beta_n(t) \notin U(\alpha(\infty), r', \epsilon')$ for all $\epsilon' < r'$.
\end{example}

In \cite{murray:2015aa} Murray does prove a weaker type of North-South dynamics on the contracting boundary.

\begin{theorem}[Corollary 4.3 \cite{murray:2015aa}] \label{thm: weak n-s} Let $X$ be a proper $\mathrm{CAT}(0)$ space and let $G$ be a group acting geometrically on $X$. If $g$ is a rank-1 isometry in $G$, i.e., a contracting element, $U$ is an open neighborhood of $g^\infty$ and $K$ is a compact set in $\partial_c X \smallsetminus g^{-\infty}$ then for sufficiently large $n$, $g^n(K) \subset U$. 
\end{theorem}

Murray also uses dynamical methods to prove a classical result that is known for the action of a hyperbolic group on its boundary. 

\begin{theorem}[Theorem 4.1  \cite{murray:2015aa}] Let $G$ be a group acting geometrically on a proper $\mathrm{CAT}(0)$ space. Either $G$ is virtually $\Z$ or the $G$ orbit of every point in the contracting boundary is dense.
\end{theorem}

In Section \ref{subsec:Morse boundary} we will introduce a generalization of the contracting boundary to any proper geodesic space. It is unknown if any of these dynamical results hold in the more general settings. 


\subsection{Morse boundary} \label{subsec:Morse boundary}
The Morse boundary, introduced by Cordes in \cite{cordes:2016ad}, generalizes the contracting boundary to the setting of proper geodesic spaces. This boundary retains many of the nice properties of the contracting boundary including quasi-isometry invariance and visibility. In the case of a proper $\mathrm{CAT}(0)$ space it is the contracting boundary, and in the case of a proper hyperbolic space it is the Gromov boundary. The generality in which this boundary is defined means it is a quasi-isometry invariant for every finitely generated group.

We will see later in Section \ref{sec:Metric Morse} that there is a more general definition of the Morse boundary of a not-necessarily-proper geodesic space. This definition will use the Gromov product and will carry more structure, but it is often helpful to reduce to the case when the boundary can be defined by geodesics instead of sequences (as we will see in Section \ref{subsec:boundary cc}) for proper geodesic spaces. 

Let $X$ be a proper geodesic space and fix a basepoint $e \in X$. The \emph{Morse boundary} of $X$, $\partial_M X$, is the set of all Morse geodesic rays in $X$ (with basepoint $e$) up to asymptotic equivalence. To topologize the boundary, first fix a Morse gauge $N$ and consider the subset of the Morse boundary that consists of all rays in $X$ with Morse gauge at most $N$:  \begin{equation*} \partial_M^N X_e= \{\alpha(\infty) \mid \exists \beta \in \alpha(\infty) \text{ that is an $N$--Morse geodesic ray with } \beta(0)=e\}. \end{equation*} Unlike in the case of a $\mathrm{CAT}(0)$ space, the visual topology of the boundary of a proper geodesic may not even be well-defined. So, instead, we take a page from the definition of the Gromov topology on the boundary of a hyperbolic space. This first step is a lemma in \cite{cordes:2016ad} that says that two $N$-Morse geodesic rays with the same basepoint which fellow travel stay uniformly close and that uniform bound, $\delta_N$, only depends on $N$. Thus we can topologize this set in a similar manner as one does for the Gromov boundary of hyperbolic spaces: the topology is defined by a system of neighborhoods, $\{V_n(\alpha) \mid n \in \N \}$, at a point $\alpha$ in $\partial_M^N X_e$. The sets $V_n( \alpha)$ are defined to be the set of geodesic rays $\gamma$ with basepoint $e$ and $d(\alpha(t), \gamma(t))< \delta_N$ for all $t<n$. That is, two $N$-Morse rays are close in $\partial_M^N X_e$ if they stay closer than $\delta_N$ for a long time.

Let $\mathcal M$ be the set of all Morse gauges. We put a partial ordering on $\mathcal M$ so that  for two Morse gauges $N, N' \in \mathcal M$, we say $N \leq N'$ if and only if $N(\lambda,\epsilon) \leq N'(\lambda,\epsilon)$ for all $\lambda,\epsilon \in \N$. We define the Morse boundary of $X$ to be
 \begin{equation*} \partial_M X_e=\varinjlim_\mathcal{M} \partial^N_M X_e \end{equation*} with the induced direct limit topology, i.e., a set $U$ is open in $\partial_M X_e$ if and only if $U \cap \partial^N_M X_e$ is open for all $N$.   
 
The Morse boundary retains almost all of the properties of the contracting boundary. The one exception is that it is open whether or not the Morse boundary is $\sigma$-compact because a priori the direct limit is over an uncountable set.
 
\begin{theorem}[\cite{cordes:2016ad}] \label{thm:Morse main} Given a proper geodesic space $X$, the Morse boundary, $\partial_M X=\varinjlim \partial_M^N X_e$, equipped with the direct limit topology, is 
\begin{enumerate}
\item a visibility space, i.e., any two points in the Morse boundary can be joined by a bi-infinite Morse geodesic;
\item independent of choice of basepoint;
\item a quasi-isometry invariant; and
\item homeomorphic to the Gromov boundary if $X$ is hyperbolic and to the contracting boundary if $X$ is $\mathrm{CAT}(0)$.
\end{enumerate}
\end{theorem}

One useful property of the Morse boundary is that compact subsets consist of uniformly Morse geodesics \cite{murray:2015aa, cordes:2016ab}. Futhermore, as in the case of the contracting boundary, a group has a compact Morse boundary if any only if it is hyperbolic \cite{cordes:2016ab}.

\section{(metric) Morse boundary and stability} \label{sec:Metric Morse}

An alternative approach to understanding ``hyperbolic directions" in a metric space is to understand ``hyperbolic" or quasi-convex subgroups/subspaces. In the case of hyperbolic groups, quasi-convex subgroups are finitely generated and undistorted. Furthermore these properties are preserved under quasi-isometry. In a general group, though, quasi-convexity depends on a choice of generating set and is not preserved by quasi-isometry. Thus in an effort to preserve these qualities, we look at a stronger notion of quasi-convexity:


\begin{definition} We say a quasi-convex subspace $Y$ of a geodesic metric space $X$ is \emph{$N$-stable} if every pair of points in $Y$ can be connected by a geodesic which is $N$-Morse in $X$. We say that a subgroup is stable if it is stable as a subspace.
\end{definition}
\begin{remark}[Relationship with Durham--Taylor definition] It is important to note that this is a generalization of the original definition of stability given by Durham--Taylor in \cite{durham:2015aa}. The definition above detects the same collection of stable subsets up to quasi-isometry, and the two definitions coincide for subgroups of finitely generated groups \cite[Lemma 3.8]{cordes:2016aa}.
\end{remark}

Durham--Taylor prove that the collection of stable subgroups of mapping class groups are precisely those which are convex-cocompact in the sense of Farb--Mosher \cite{farb:2002aa,durham:2015aa}. These subgroups are well studied and have important connections to the geometry of Teichm\"uller space, the curve complex and surface group extensions. In the setting of right-angled Artin groups, Koberda, Manghas, and Taylor classify all the stable subgroups. They prove that these subgroups are all free \cite{koberda:2014aa}. For more on the history and current situation of stable subgroups see Section \ref{sec:stable subgroups}.
\medskip

We will see in this section that the notions of the Morse boundary and of stability can be united. We will do this by viewing any geodesic metric space as the union of stable subsets which are indexed by Morse gauges $N$ and hyperbolic (with hyperbolicity constant depending only on $N$). We will define these subspaces in Section \ref{sec:stable strata}, but before this we will recall some definitions.

\subsection{sequential boundary, capacity dimension \& asymptotic dimension} \label{sec:seq, capdim, asdim}

\begin{definition}
Let $X$ be a metric space and let $x,y,z\in X$. The \emph{Gromov product} of $x$ and $y$ with respect to $z$ is defined as $$ (x\cdot y)_z = \frac{1}{2}\left( d(z,x)+d(z,y)-d(x,y)\right).$$
Let $(x_n)$ be a sequence in $X$. We say $(x_n)$ converges at infinity if $(x_i\cdot x_j)_e \to \infty$ as $i,j \to \infty$. Two convergent sequences $(x_n),(y_m)$ are said to be \emph{equivalent} if $(x_i \cdot y_j) \to \infty$ as $i,j \to \infty$. We denote the equivalence class of $(x_n)$ by $\lim x_n$.

The \emph{sequential boundary} of $X$, denoted $\partial X$, is defined to be the set of convergent sequences considered up to equivalence.
\end{definition}

\begin{definition}[4-point definition of hyperbolicity; Definition $1.20$ \cite{bridson:1999fj}]\label{defn:hyp}  Let $X$ be a (not necessarily geodesic) metric space. We say $X$ is $\delta$\emph{--hyperbolic} if for all $w,x,y,z$ we have
$$ (x\cdot y)_w \geq \min\set{(x\cdot z)_w,(z\cdot y)_w} - \delta.$$
\end{definition}

If  $X$ is $\delta$--hyperbolic, we may extend the Gromov product to $\partial  X$ in the following way:
$$ (x\cdot y)_e = \sup \left(\liminf_{m,n\to\infty}\set{(x_n\cdot y_m)_e}\right).$$
where $x,y \in \partial X$ and the supremum is taken over all sequences $(x_i)$ and $(y_j)$ in $X$ such that $x= \lim x_i$ and $y=\lim y_j$.

Recall that a metric $d$ on $\partial  X$ is said to be \emph{visual} (with parameter $\varepsilon>0$) if there exist $k_1,k_2>0$ such that $k_1\exp(-\varepsilon(x\cdot y)_e) \leq d(x,y)\leq k_2\exp(-\varepsilon(x\cdot y)_e)$, for all $x,y\in\partial X$.

Let $x, y \in \partial  X$. As a shorthand we define $\rho_\epsilon(x, y) := \exp\left(-\varepsilon (x\cdot y\right)_e)$.

\begin{theorem}\label{thm:GrProdmetric} \textrm{\cite[Section $7.3$]{Ghys:1990aa}} Let $X$ be a $\delta$--hyperbolic space. If $\varepsilon'=\exp(2\delta\varepsilon)-1\leq \sqrt{2}-1$ then we can construct a visual metric $d$ on $\partial  X$ such that
\[
 (1-2\varepsilon')\rho_\varepsilon(x, x')\leq d(x,x') \leq \rho_\varepsilon(x,x').
\]
\end{theorem}

Visual metrics on a hyperbolic space are all quasi-symmetric. A quasi-symmetry is a map that is a generalization of bi-Lipschitz map: instead of controlling how much the diameter of a set can change, a quasi-symmetry preserves only the relative sizes of sets.

\begin{definition}\label{defn:quasisym} A homeomorphism $f\colon(X,d)\to (Y,d')$ is said to be quasi--symmetric if there exists a homeomorphism $\eta:\R\to\R$ such that for all distinct $x,y,z\in X$, 
\[
 \frac{\displaystyle d'(f(x),f(y))}{\displaystyle d'(f(x),f(z))} \leq \eta \left( \frac{\displaystyle d(x,y)}{\displaystyle d(x,z)}\right).
\]
\end{definition}

One natural quasi-isometry invariant assigned the boundary of a hyperbolic space is the \emph{capacity dimension}. This invariant was introduced by Buyalo in \cite{Buyalo:2005aa} and it is sometimes known as linearly-controlled dimension.

Let $\mathcal{U}$ be an open covering of a metric space $X$. Given $x \in X$, we let \[L(\mathcal{U},x)=\sup\setcon{d(x, X\char`\\U)}{U \in \mathcal{U}}\] be the \emph{Lebesgue number of $\mathcal{U}$ at $x$ } and $L(\mathcal{U}) = \inf_{x \in X} L(\mathcal{U},x)$ the Lebesgue number of $\mathcal{U}$. The \emph{multiplicity of $\mathcal{U}$}, $m(\mathcal{U})$, is the maximal number of members of $\mathcal{U}$ with non-empty intersection.

\begin{definition}[\cite{Buyalo:2007ab}] The \emph{capacity dimension} of a metric space $X$, $\cdim{X}$, is the minimal integer $m$ with the following property:

There exists some $\delta \in (0,1)$ such that for every sufficiently small $r>0$ there is an open covering $\mathcal{U}$ of $X$ by sets of diameter at most $r$ with $L(\mathcal{U}) \geq \delta r$ and $m(\mathcal{U})\leq m+1$.
\end{definition}

The capacity dimension is similar to the covering dimension in that it is an infimum over open covers, but the capacity dimension necessitates metric information: given an open cover $\mathcal{U}$ the capacity dimension requires a linear relationship between the $\sup_{U \in \mathcal{U}}\{\mathrm{diam}(U)\}$ and $L(\mathcal{U})$.  For more information see \cite{Buyalo:2007ab}.

In \cite[Corollary 4.2]{Buyalo:2005aa}, Buyalo shows that the capacity dimension of a metric space is a quasi-symmetry invariant and since quasi-isometries of hyperbolic spaces induce quasi-symmetries on the boundaries, this shows that the capacity dimension of the boundary is an invariant of a hyperbolic space.

Another quasi-isometry invariant dimension one can assign to a metric space is the asymptotic dimension. This notion was introduced by Gromov in \cite{gromov:1993aa} and is a coarse version of the topological dimension.
\begin{definition} A metric space $X$ has  \emph{asymptotic dimension at most} $n$ (${asdim}(X)\leq n$), if for every $R>0$ there exists a cover of $X$ by uniformly bounded sets such that every metric $R$--ball in $X$ intersects at most $n+1$ elements of the cover. We say $X$ has \emph{asymptotic dimension $n$} if ${asdim}(X) \leq n$ but ${asdim}(X) \nleq n-1$.
\end{definition}
 The celebrated theorem of Yu showed that groups with finite asymptotic dimension satisfy both the coarse Baum--Connes and the Novikov conjectures \cite{yu:1998aa}. Many classes of groups have been shown to have finite asymptotic dimension including: hyperbolic groups \cite{Roe:2005aa}, relatively hyperbolic groups whose parabolic subgroups are of finite dimension \cite{osin:2005aa}, mapping class groups \cite{bestvina:2015ab}, cubulated groups \cite{wright:2012aa}. But exact bounds are often hard to calculate. 

Buyalo and Lebedeva used the capacity dimension to prove a conjecture of Gromov: the asymptotic dimension of any hyperbolic group is the topological dimension of its boundary plus one \cite{Buyalo:2007aa}. 

\subsection{stable strata} \label{sec:stable strata}

\begin{definition}[stable stratum] Let $X$ be a geodesic metric space and  $e \in X$. We define $X^{(N)}_e$ to be the set of all points in $X$ which can be joined  to $e$ by an $N$-Morse geodesic. 
\end{definition}

What do these stable strata look like? First, we can easily see that given $x \in X$ there exists  $N$ such that $x \in X^{(N)}_e$ by choosing $N=N(\lambda, \epsilon)=\lambda d(e,x) + \epsilon$. So the collection of $X^{(N)}_e$ will cover $X$. If $X$ is hyperbolic, since all geodesics are $N$-Morse for some $N$, we have that $X=X^{(N)}_e$ for that $N$. On the other hand, if $X$ is a space with no infinite Morse geodesics (e.g., groups satisfying non-trivial laws or products of unbounded spaces) then the $X^{(N)}_e$ are just sets of bounded diameter. For groups with mixed geometries this is a hard question to answer. We will see later in the survey (Sections \ref{sec:calcs in fg grps} and \ref{sec:rel hyp grps})  that (up to quasi-isometry) we can begin to understand these subspaces for some groups.

One thing we do know is that these subspaces are hyperbolic. A standard argument shows that if $x,y \in X^{(N)}_e$ then the geodesic triangle in $X$ formed by $x,y,e$ is $4N(3,0)$-slim \cite[Lemma 2.2]{cordes:2016ad} proving these spaces are hyperbolic. Note, though, that the geodesic $[x,y]$ is not necessarily contained in $X^{(N)}_e$.  So since the $ X^{(N)}_e$ are not necessarily geodesic, we use the 4-point definition of hyperbolicity and conclude they are hyperbolic.


Using the fact that these triangles are slim, one can also show that if $x,y \in X^{(N)}_e$ then $[x,y]$ is $N'$-Morse where $N'$ depends only on $N$ \cite[Lemma 2.3]{cordes:2016ad}. This shows that the $X^{(N)}_e$ are not only quasi-convex, but they are $N'$-stable!

Furthermore, as each $X^{(N)}_e$ is hyperbolic we may consider its Gromov boundary, $\partial_s X^{(N)}_e$, and the associated visual metric $d_{(N)}$. (See \cite[Section 2.2]{Buyalo:2007ab} for a careful treatment of the sequential boundary of a hyperbolic space which is not necessarily geodesic.) This metric is unique up to quasi-symmetry. 

Natural maps between strata have nice properties: the natural inclusion $X^{(N)}_e\subseteq X^{(N')}_e$ induces a map $\partial_s X^{(N)}_e\to \partial_s X^{(N')}_e$ which is a quasi-symmetry onto its image. Additionally, if we have a quasi-isometry $q\colon X\to Y$, then for every $N$ there exists an $N'$ such that $q(X^{(N)}_e)\subseteq Y^{(N')}_{q(e)}$. This induces an embedding $\partial q \colon \partial_s X^{(N)}_e\to \partial_s Y^{(N')}_{q(e)}$ which is a quasi-symmetry onto its image. We will see in Section \ref{sec:stable equivalence invariants} that these properties will be useful in defining new quasi-isometry invariants.

Finally, if $X$ is also a proper metric space, then for each $N$ there is a homeomorphism between $\partial_M^N X$ (as defined in Section \ref{subsec:Morse boundary} with geodesics) and $\partial_s X^{(N)}_e$ and thus the Morse boundary is homeomorphic to direct limit over $\partial_s X^{(N)}_e$ as topological spaces.

In summary:
\begin{theorem}[\cite{cordes:2016aa}]\label{bthm:key} Let $X,Y$ be geodesic metric spaces and let $e\in X$. The family of subsets $X^{(N)}_e$ of $X$ indexed by functions $N\colon \R_{\geq1} \times \R_{\geq0}\to\R_{\geq0}$ enjoys the following properties:
\begin{enumerate}[label=\textup{\Roman*}]
\item (covering) $X=\bigcup_N X^{(N)}_e$.\label{key:covering}
\item (partial order) If $N\leq N'$, then $X^{(N)}_e\subseteq X^{(N')}_e$. \label{key:PO}
\item (hyperbolicity) Each $X^{(N)}_e$ is hyperbolic in the sense of Definition $\ref{defn:hyp}$.\label{key:hyperbolic}
\item (stability) Each $X^{(N)}_e$ is $N'$-stable, where $N'$ depends only on $N$.\label{key:stable}
\item (universality) Every stable subset of $X$ is a quasi-convex subset of some $X^{(N)}_e$.\label{key:universal}
\item (boundary) The sequential boundary $\partial_s X^{(N)}_e$ can be equipped with a visual metric which is unique up to quasi-symmetry. An inclusion $X^{(N)}_e\subseteq X^{(N')}_e$ induces a map $\partial_s X^{(N)}_e\to \partial_s X^{(N')}_e$ which is a quasi-symmetry onto its image.\label{key:boundary}
\item (generalizing the Gromov boundary) If $X$ is hyperbolic, then $X=X^{(N)}_e$ for all $N$ sufficiently large, and $\partial_s X^{(N)}_e$ is quasi-symmetric to the Gromov boundary of $X$. \label{key:Gromov}
\item (generalizing the Morse boundary) If $X$ is proper, then its Morse boundary is equal to the direct limit of the $\partial_s X^{(N)}_e$ as topological spaces.\label{key:Mboundary}
\item (behavior under quasi-isometry) If $q \colon X\to Y$ is a quasi-isometry then for every $N$ there exists an $N'$ such that $q(X^{(N)}_e)\subseteq Y^{(N')}_{q(e)}$ and there is an induced embedding $\partial q \colon \partial_s X^{(N)}_e\to \partial_s Y^{(N')}_{q(e)}$ which is a quasi-symmetry onto its image.\label{key:quasiisom}
\end{enumerate}
\end{theorem}


There is a stronger version of \hyperref[bthm:key]{\ref{bthm:key}.\ref{key:quasiisom}}, but we will first we require some additional terminology.

Let $X,Y$ be geodesic metric spaces, let $x\in X$ and $y\in Y$. We say $X$ is \emph{stably subsumed} by $Y$ (denoted $X\hookrightarrow_s Y$) if, for every $N$ there exists some $N'$ and a quasi-isometric embedding $X^{(N)}_x\to Y^{(N')}_y$. By \hyperref[bthm:key]{\ref{bthm:key}.\ref{key:quasiisom}}, this property is independent of the choice of $x,y$. We say $X$ and $Y$ are \emph{stably equivalent} (denoted $X\sim_s Y$) if they stably subsume each other. It is easy to see that two spaces are stably equivalent if and only if they have the same collection of stable subsets up to quasi-isometry.

Given a geodesic metric space $X$, we will consider the collection of boundaries $\left( \partial_s X^{(N)}_e\right)$ equipped with visual metrics as the \emph{metric Morse boundary} of $X$. 

We say that one collection of spaces $\left( A_i\right)_{i\in I}$ is \emph{quasi-symmetrically subsumed} by another $\left(B_j\right)_{j\in J}$ (denoted $(A_i)\hookrightarrow_{qs}(B_j)$) if, for every $i$ there exists a $j$ and an embedding $A_i \to B_j$ which is a quasi--symmetry onto its image. Two collections are \emph{quasi-symmetrically equivalent} (denoted $(A_i)\sim_{qs}(B_j)$) if $(A_i)\hookrightarrow_{qs}(B_j)$ and $(A_i)\hookrightarrow_{qs}(A_i)$).

\begin{ctheorem} \label{bthm:keyquasiisom2} Let $X,Y$ be geodesic metric spaces, let $x\in X$ and $y\in Y$. Then $X\hookrightarrow_s Y$ if and only if $\left( \partial X^{(N)}_x\right)\hookrightarrow_{qs} \left( \partial Y^{(N)}_y\right)$.
\end{ctheorem}

\begin{corollary}\label{bcor:stableinvarience} Quasi-isometric geodesic metric spaces are stably equivalent and have quasi-symmetrically equivalent metric Morse boundaries. In particular, the metric Morse boundary is invariant under change of basepoint.
\end{corollary}

Stable equivalence is a much weaker notion than quasi-isometry. All spaces with no infinite Morse geodesic rays will be stably equivalent to a point! On the other hand, Cordes--Hume show that the mapping class group and Teichm\"uller space with the Teichm\"uller metric are stably equivalent (Theorem \ref{thm:MCG}), when it is well known that they are not quasi-isometric.

\subsection{stable equivalence invariants} \label{sec:stable equivalence invariants} 
We can now define two stable equivalence invariants using the notions discussed in Section \ref{sec:seq, capdim, asdim}; in particular, the capacity dimension and the asymptotic dimension. The definitions will rely on Theorem \ref{bthm:key}.

\begin{definition}[stable asymptotic dimension] The \emph{stable asymptotic dimension} of $X$ ($\stabdim (X)$) is the supremal asymptotic dimension of a stable subset of $X$ \end{definition}

We can see that by universality it is possible to just consider the maximal asymptotic dimension of the $X^{(N)}_e$. One obvious but useful bound is that the stable asymptotic dimension is bounded from above by the asymptotic dimension.

\begin{definition}[Morse capacity dimension] The \emph{Morse capacity dimension} of $X$ ($\Mcd(X)$) is the supremal capacity dimension of spaces in the metric Morse boundary. We say that the empty set has capacity dimension $-1$.\end{definition}

By Theorem \ref{bthm:key}.\ref{key:boundary} we know that the inclusion map $X^{(N)}_e \into X^{(N')}_e$ induces a map $\partial X^{(N)}_e \into \partial X^{(N')}_e$ which is a quasi-symmetry onto its image. So this definition is well defined.

\medskip

It follows from Theorem \ref{bthm:key} that these notions are invariant under change of basepoint and a stable equivalence invariant and thus a quasi-isometry invariant. It also follows that  the stable dimension of a hyperbolic space is precisely its asymptotic dimension and the Morse capacity dimension of a hyperbolic space is the capacity dimension of its boundary equipped with some visual metric.

\begin{remark}[conformal dimension] Since the conformal dimension (introduced by Pansu in \cite{Pansu:1989aa}) is also a quasi-symmetry invariant, one could also define a conformal dimension of the Morse boundary in the same manner and with the properties listed above, but it is often harder to compute.
\end{remark}

By using bounds bounds proved in the hyperbolic setting, we get a nice relationship between these two dimensions \cite[Theorem $1.1$]{Buyalo:2005aa}, \cite[Proposition $3.6$]{Mackay:2013aa}.

 \begin{proposition}\label{bcor:comparedims} Let $X$ be a geodesic metric space. Then
\[
 \stabdim (X)-1\leq \Mcd(X) \leq \stabdim (X).
\] 
\end{proposition}

\subsection{calculations in finitely generated groups}\label{sec:calcs in fg grps}

We now present some calculations in finitely generated groups where Morse geodesics have been characterized in some way: mapping class groups, right-angled Artin groups, and $C'(1/6)$ graphical small cancellation groups.  As we will see, a recurring way to calculate an upper bound on the stable asymptotic dimension is to show that each stable stratum embeds into a hyperbolic space that is easier to understand.

 \subsubsection{mapping class group and Teichm\"uller space} 

Assume that $\Sigma$ is an orientable surface of finite type. Let $\teich$ be the Teichm\"uller space of $\Sigma$ with the Teichm\"uller metric. The following result is a generalization of a result in \cite{cordes:2016ad} which shows that $\partial_M \Mod$ is homeomorphic to $\partial_M \teich$.
 
 \begin{theorem}[\cite{cordes:2016aa}]\label{thm:MCG} $\mathrm{MCG}(\Sigma)$ and $\teich$ are stably equivalent thus $\stabdim(\mathrm{MCG}(\Sigma))=\stabdim(\teich)$. Furthermore $\stabdim(\mathrm{MCG}(\Sigma))$ and $\Mcd{\Mod}$ are bounded above linearly in the complexity of $\Sigma$. \end{theorem}
 
 An upper bound on the stable asymptotic dimension of mapping class groups can be obtained via the bounds on the asymptotic dimension for mapping class groups obtained by Bestvina--Bromberg--Fujiwara in \cite{bestvina:2015ab} or by Behrstock--Hagen--Sisto in \cite{behrstock:2015aa} which are exponential or quadratic in the complexity of the surface respectively. To show there is a linear bound on the stable asymptotic dimension, Cordes--Hume show each stable subset of $\Mod$ quasi-isometrically embeds into the curve graph. They then use the bound found by Bestvina--Bromberg on the asymptotic dimension of the curve graph \cite{bestvina:2015aa}.

Leininger and Schleimer prove that for every $n$ there is a surface $\Sigma$ such that the Teichm\"{u}ller space $\teich$ contains a stable subset quasi--isometric to $\mathbb{H}^n$  \cite{leininger:2014aa}. This fact gives a lower bound on the stable dimension, which is at best logarithmic in the complexity, but also shows that Teichm\"uller spaces can have arbitrarily high stable asymptotic dimension. Since $\teich$ is stably equivalent to $\Mod$, we see that $\Mod$ contains a stable subset quasi-isometric to $\mathbb{H}^n$. The only known explicit examples of convex cocompact subgroups of mapping class groups are virtually free groups \cite{Dowdall:2014ab, Kent:2009aa, koberda:2014aa, Min:2011aa}. Although the results of Leininger--Schleimer do not provide any non-virtually-free convex cocompact subgroups, the fact that $\stabdim(\mathrm{Mod} (\Sigma))>1$ for some surfaces shows that there is no purely geometric obstruction to the existence of a non-free convex cocompact subgroup of $\Mod$.

\subsubsection{right-angled Artin groups}

Koberda, Manghas, and Taylor classify the stable subgroups of all right-angled Artin groups \cite{koberda:2014aa}. They show that these subgroups are always free. The next theorem is the natural analogue for stable subspaces. 

\begin{theorem}[\cite{cordes:2016aa}] \label{thm:RAAG} Let $X$ be a Cayley graph of a right--angled Artin group. Every stable subset of $X$ is quasi-isometric to a proper tree. In particular, $X$ is stably equivalent to a line if the group is Abelian of rank $1$, a point if it is Abelian of rank $\neq 1$ and a regular trivalent tree otherwise.
\end{theorem}

The proof of this theorem follows in a very similar manner to the proof of Theorem \ref{thm:MCG}. In this case each stable subset of $X$ embeds into the contact graph (defined by Hagen \cite{hagen:2014aa}), which is a quasi-tree. As a result, each $X^{(N)}_e$ is quasi-isometric to a proper tree. The proof is completed by calling on the universality condition in Theorem \ref{bthm:key}.  By universality of stable subsets (Theorem \hyperref[bthm:key]{\ref{bthm:key}.\ref{key:universal}}), we know that any stable subgroup of a right-angled Artin group is quasi-isometric to a proper simplicial tree. Thus since groups which are quasi-isometric to trees are virtually free  \cite[Corollary 7.19]{Ghys:1990aa}, they recover a result of \cite{koberda:2014aa}: stable subgroups of right--angled Artin groups are free. 



\subsubsection{graphical small cancellation}

We move to the realm of graphical small cancellation groups. Graphical small cancellation theory is a generalization introduced by Gromov \cite{gromov:2003aa} in order to construct groups whose Cayley graphs contain certain prescribed subgraphs, in particular one can construct ``Gromov monster'' groups, those with a Cayley graph which coarsely contains expanders \cite{Arzhantseva:2008aa, Osajda:2014aa}. These monster groups cannot be coarsely embedded into a Hilbert space, and they are the only known counterexamples to the Baum--Connes conjecture with coefficients \cite{Higson:2002aa}.

\begin{theorem}[\cite{cordes:2016aa}] \label{bthm:smallcanc} Let $X$ be the Cayley graph of a graphical $C'(1/6)$ small cancellation group. Then $\stabdim (X)\leq 2$ and $\Mcd(X)\leq 1$.
\end{theorem}

Note that this is optimal as fundamental groups of higher genus surfaces are hyperbolic with asymptotic dimension $2$ and admit $C'(\frac{1}{6})$ graphical small cancellation presentations.

Again, we work with the stable strata. Each of the strata embeds quasi-isometrically into a finitely presented classical $C'(1/6)$ small cancellation group. These groups are hyperbolic with asymptotic dimension at most $2$ and the capacity dimension of their Gromov boundaries is at most $1$.




\subsection{relatively hyperbolic groups} \label{sec:rel hyp grps}
 Osin in \cite{osin:2005aa} shows that relatively hyperbolic groups inherit finite asymptotic dimension from their maximal parabolic subgroups. This is also true for the stable asymptotic dimension:

\begin{theorem}[\cite{cordes:2016aa}] \label{bthm:relhyp} Let $G$ be a finitely generated group which is hyperbolic relative to $\mathcal H$. Then $\stabdim (G)<\infty$ if and only if $\stabdim (H)<\infty$ for all $H\in\mathcal{H}$.
\end{theorem}

In \cite{cordes:2016ac} Cordes--Hume focus on relatively hyperbolic groups. In this paper they suggest an approach to answering the following question which appears in \cite{Behrstock:2009aa}: How may we distinguish non-quasi-isometric relatively hyperbolic groups with non-relatively hyperbolic peripheral subgroups when their peripheral subgroups are quasi-isometric?

Using small cancellation theory over free products Cordes--Hume construct quasi-isometrically distinct one-ended relatively hyperbolic groups which are all hyperbolic relative to the same collection of groups. These groups are distinguished using a notion similar to stable asymptotic dimension; rather than $X^{(N)}_e$, we use stable subsets which ``avoid" the left cosets of the peripheral subgroups.

\begin{theorem}[\cite{cordes:2016ac}] Let $\mathcal H$ be a finite collection of finitely generated groups each of which has finite stable dimension or are non-relatively hyperbolic. There is an infinite family of $1$-ended groups $(G_n)_{n\in\N}$, which are non-quasi-isometric, where each $G_n$ is hyperbolic relative to $\mathcal H$.
\end{theorem}

\section{stable subgroups} \label{sec:stable subgroups}

We will start our discourse on stable subgroups with some with some motivation from Kleinian groups and mapping class groups.

A non-elementary discrete (Kleinian) subgroup $\Gamma < \PSL_2(\C)$ determines a minimal $\Gamma$-invariant closed subspace $\Lambda(\Gamma)$ of the Riemann sphere called its \emph{limit set} and taking the convex hull of $\Lambda(\Gamma)$ determines a convex subspace of $\mathbb H^3$ with a $\Gamma$-action.  A Kleinian group $\Gamma$ is called \emph{convex cocompact} if it acts cocompactly on this convex hull or, equivalently, any $\Gamma$-orbit in $\mathbb H^3$ is quasiconvex.  

Originally defined by Farb--Mosher \cite{farb:2002aa} and later developed further by Kent--Leininger \cite{Kent:2008aa} and Hamenst\"adt \cite{Hamenstadt:2005aa}, a subgroup $H < \Mod$ is called convex cocompact if and only if any $H$-orbit in $\teich$, the Teichm\"uller space of $\Sigma$ with the Teichm\"uller metric, is quasiconvex, or $H$ acts cocompactly on the weak hull of its limit set $\Lambda(H) \subset \mathbb{P}\mathcal{MF}(\Sigma)$ in the Thurston compactification of $\teich$.  This notion is important because convex cocompact subgroups $H< \Mod$ are precisely those which determine Gromov hyperbolic surface group extensions. Furthermore, Farb--Mosher show that if there is a purely pseudo-Anosov subgroup of $\Mod$ that is not convex cocompact, then this subgroups would be a counterexample to Gromov's conjecture that every group with a finite Eilenberg--Mac Lane space and no Baumslag--Solitar subgroups is hyperbolic (see \cite{Kent:2007aa} for more information).

In both of these examples, convex cocompactness is characterized equivalently by both a quasi-convexity condition and an asymptotic boundary condition.  In \cite{durham:2015aa}, Durham--Taylor introduced stability in order to characterize convex cocompactness in $\Mod$ by a quasiconvexity condition intrinsic to the geometry of $\Mod$ and this condition naturally generalizes the above quasi-convexity characterizations of convex cocompactness to any finitely generated group. There has been much recent work to characterize stable subgroups of important groups. The theorem below is a brief summary of this work.

\begin{theorem} \label{thm: stable summary} Let the pair $(G,H)$ of a finitely generated group $G$ and a subgroup $H$ satisfy one of the following:
\begin{enumerate}
\item $G$ is hyperbolic and $H$ is quasi-convex; \label{qc in hyp}
\item G is relatively hyperbolic and $H$ is a finitely generated and quasi-isometrically embeds in the the coned off graph in the sense of \cite{Farb:1998aa}; \label{rel hy stable}
\item $G=A(\Gamma)$ is a right-angled Artin group with $\Gamma$ a finite graph which is not a join and $H$ is a finitely generated subgroup quasi-isometrically embedded in the extension graph \cite{koberda:2014aa}; \label{raag stab}
\item $G=\Mod$ and $H$ is a convex cocompact subgroup in the sense of \cite{farb:2002aa}; \label{mcg stab}
\item $G=\Out$ and $H$ is a convex cocompact subgroup in the sense of \cite{Hamenstadt:2014aa}; \label{out stab}
\item $H$ is hyperbolic and hyperbolically embedded in $G$; \label{hy emb stab}
\end{enumerate}
Then $H$ is stable in $G$. Moreover for (\ref{qc in hyp}), (\ref{raag stab}), and (\ref{mcg stab}), the reverse implication holds.
\end{theorem}

Item (\ref{qc in hyp}) follows easily by checking the definition. Item (\ref{rel hy stable}) is due to \cite{Aougab:2016aa}, for (\ref{raag stab}) this is due to \cite{koberda:2014aa}, for (\ref{mcg stab}) this is due to \cite{durham:2015aa}, for (\ref{out stab}) this is again \cite{Aougab:2016aa}, and for (\ref{hy emb stab}) this follows from \cite{sisto:2016aa}.

\subsection{boundary cocompactness} \label{subsec:boundary cc}

There is an alternative characterization to the Durham--Taylor notion of stability using the Morse boundary to define an asymptotic property for subgroups of finitely generated groups called \emph{boundary convex cocompactness} which generalizes the classical boundary characterization of convex cocompactness from Kleinian groups. This is presented in \cite{cordes:2016ab}.

Let $G$ be a finitely generated group acting by isometries on a proper geodesic metric space $X$.  Fix a basepoint $e \in X$.  In order to define this boundary characterization we first need to define a limit set.

\begin{definition}[$\Lambda(G)$] \label{defn:limit set}
The \emph{limit set} of the $G$-action on $\partial_M X$ is
$$\Lambda_e(G)= \setcon{\lambda \in \partial_M X}{\exists N \text{ and } (g_k) \subset G  \text{ such that } (g_k \cdot e)\subset X^{(N)}_e \text{ and } \lim g_k\cdot e=\lambda}$$
\end{definition}

That is, the limit set $\Lambda_e(G)$ is the set of points which can be represented by sequences of uniformly Morse $G$-orbit points;  note that $\Lambda_e(G)$ is obviously $G$-invariant. It is also invariant under change of basepoint. 

Given some $\lambda \in \Lambda_e(G)$ we know there is a sequence $(g_k \cdot e)$ so that $\lim g_k\cdot e=\lambda$. We can produce a geodesic from this sequence in a standard way: for each $k$ let $\alpha_k$ be an $N$-Morse geodesic joining $e$ and $g_k \cdot e$. It follows from Arzel\`a--Ascoli that there is a subsequence $(\alpha_{i(k)})$ which converges (uniformly on compact sets) to a geodesic ray that is $N$-Morse by \cite[Lemma 2.8]{cordes:2016ad}. This map defines a homeomorphism between the Morse boundary defined by sequences and the Morse boundary defined by geodesics \cite[Theorem 3.14]{cordes:2016aa}. From this perspective, by Theorem \ref{thm:Morse main} we know that there exists a bi-infinite Morse geodesic joining any two distinct points in the limit set. This is a starting point for defining the weak hull, but we want to define it by taking \emph{all} geodesics with distinct endpoints in $\Lambda_e(G)$. Formalizing this motivates the next definition:




\begin{definition}[asymptotic, bi-asymptotic]
Let $\gamma \colon (-\infty, \infty)\rightarrow X$ be a biinfinite geodesic in $X$ with $\gamma(0)$ a closest point to $e$ along $\gamma$.  Let $\lambda \in \partial X^{(N)}_e$.  We say $\gamma$ is \emph{forward asymptotic} to $\lambda$ if for any $N$-Morse geodesic ray $\gamma_{\lambda}:[0, \infty)\rightarrow X$ with $\gamma_{\lambda}(0)=e$, there exists $K>0$ such that
$$d_{Haus}(\gamma([0, \infty)), \gamma_{\lambda}([0, \infty)))<K.$$
We define \emph{backwards asymptotic} similarly.  If $\gamma$ is forwards, backwards asymptotic to $\lambda, \lambda'$, respectively, then we say $\gamma$ is \emph{bi-asymptotic} to $(\lambda, \lambda')$.
\end{definition}

Now for the definition of weak hull.

\begin{definition}[weak hull]
The \emph{weak hull} of $G$ in $X$ based at $e \in X$, denoted $\mathfrak{H}_e(G)$, is the collection 
 of all biinfinite rays $\gamma$ which are bi-asymptotic to $(\lambda, \lambda')$ for some $\lambda \neq \lambda' \in \Lambda_e(G)$.
\end{definition}

An important fact about the weak hull is that of  $|\Lambda_e(G)| \geq 2$, then $\mathfrak H_e(G)$ is nonempty and $G$-invariant.

We are now ready to define boundary convex cocompactness:

\begin{definition}[boundary convex cocompactness] \label{defn:cc subspace}
We say that $G$ acts \emph{boundary convex cocompactly} on $X$ if the following conditions hold:
\begin{enumerate}
\item $G$ acts properly on $X$;
\item For some (any) $e \in X$, $\Lambda_e(G)$ is nonempty and compact; \label{compact}
\item For some (any) $e \in X$, the action of $G$ on $\mathfrak H_e(G)$ is cobounded. 
\end{enumerate} 
\end{definition}

\begin{definition}[Boundary convex cocompactness for subgroups] \label{defn:cc subgroup}
Let $G$ be a finitely generated group.  We say $H<G$ is \emph{boundary convex cocompact} if $H$ acts boundary convex cocompactly on any Cayley graph of $G$ with respect to a finite generating set. 
\end{definition}

\begin{theorem}[\cite{cordes:2016ab}]\label{thm:main}
Let $G$ be a finitely generated group.  Then $H<G$ is boundary convex cocompact if and only if $H$ is stable in $G$.
\end{theorem}

Theorem \ref{thm:main} is an immediate consequence of the following stronger statement:

\begin{theorem}[\cite{cordes:2016ab}] \label{thm:cc and stab}
Let $G$ be a finitely generated group acting by isometries on a proper geodesic metric space $X$.  Then the action of $G$ is boundary convex cocompact if and only if some (any) orbit of $G$ in $X$ is a stable embedding.

In both cases, $G$ is hyperbolic and any orbit map $orb_e\colon G \rightarrow X$ extends continuously and $G$-equivariantly to an embedding of $\partial_{Gr} G$ which is a homeomorphism onto its image $\Lambda_e(G) \subset \partial_M X_e$.
\end{theorem}

We note that Theorem \ref{thm:main} and \cite[Proposition 3.2]{durham:2015aa} imply that boundary convex cocompactness is invariant under quasi-isometric embeddings.

\begin{remark}[The compactness assumption on $\Lambda_e(G)$ is essential for Theorem \ref{thm:cc and stab}]

Consider the group $G=\Z^2 * \Z * \Z = \fgen{a,b} * \fgen{c} * \fgen{d}$ acting on its Cayley graph. Consider the subgroup $H=\fgen{a,b,c}$.  Since the $H$ is isometrically embedded and convex in $G$, it follows that $\partial_M H_e\cong \Lambda_e(H) \subset \partial_M G_e$ and $\mathfrak H_e(H) = H$ for any $e \in G$, whereas $H$ is not hyperbolic and thus not stable in $G$. In fact, the compactness assumption ensures that the weak hull will be a subspace of some $X^{(N)}_e$, i.e., hyperbolic.
\end{remark}

The following is an interesting question:

\begin{question}
If $|\Lambda(G)|\neq \emptyset$, then must we have in fact $|\Lambda(G)|\geq 2$?
\end{question}

Recall from Section \ref{sec:dynamics} that in the case of $\mathrm{CAT}(0)$ spaces, this question has been affirmatively answered \cite[Lemma 4.9]{murray:2015aa}.

\subsubsection{height, width, bounded packing}

Antol\'in, Mj, Sisto and Taylor in \cite{Antolin:2016aa} use the boundary cocompactness characterization to extend some well-known intersection properties of quasi-convex subgroups of hyperbolic or relatively hyperbolic groups \cite{Gitik:1998aa, Hruska:2009aa} to the context of stable subgroups of finitely generated groups:

\begin{theorem} Let $H_1, \ldots H_l$ be stable subgroups of a finitely generated group. Then the collection $H = \{H_1,\ldots ,H_l\}$ satisfies the following:
\begin{enumerate}
\item H has finite height.
\item H has finite width.
\item H has bounded packing.
\end{enumerate}
\end{theorem}

In particular they show any group-subgroup pair $(G,H)$ satisfying one of the conditions in Theorem \ref{thm: stable summary} then $H$ has finite height, finite width, and bounded packing.

\section{a metrizable topology on the Morse boundary} \label{sec: Metrizing}

Cashen and Mackay have introduced a topology on the Morse boundary of a proper geodesic space that is metrizable \cite{Cashen:2017aa}. They use a generalized notion of contracting geodesics which follows that of Arzhantseva, Cashen, Gruber, and Hume \cite{Arzhantseva:2016aa}. We present the definition here:

\begin{definition} We call a function $\rho$ \emph{sublinear} if it is non-decreasing, eventually non-negative, and $\lim_{r \to \infty} \rho(r)/r =0$. \end{definition}

\begin{definition} Let $X$ be a proper geodesic metric space. Let $\gamma \colon [0, \infty) \to X$ be a geodesic ray, and let $\pi_\gamma$ be the closest point projection to $\gamma$. Then, for a sublinear function $\rho$, we say that $\gamma$ is \emph{$\rho$-contracting} if for all $x$ and $y$ in $X$: \[ d(x, y) \leq d(x, \gamma) \implies \mathrm{diam}(\pi_\gamma(x) \cup \pi_\gamma(y)) \leq \rho(d(x, \gamma)). \] \end{definition}

We see that Definition \ref{def: str contracting} is simply this definition if we ask that $\rho$ is the constant function $D$. We revisit the example in Remark \ref{rem: contracting neq morse}: the space $X$ was a ray $Y$ with set of intervals $\{I_i\}$ of length $i$ with an interval of length $i^2$ attached to the endpoints of $I_i$. We noted that this ray was not strongly contracting. It is not hard to see that in this more general definition, that it is $\sqrt{2x}$-contracting. We showed in Remark \ref{rem: contracting neq morse} that this ray was Morse. This fact is no coincidence; in proper geodesic spaces $\rho$-contracting is equivalent to being Morse \cite{Arzhantseva:2016aa}.

Cashen and Mackay introduce a new topology on the contracting boundary that is finer than the ``subspace topology" defined with the Gromov product but less fine than the direct limit topology. They call this topology the \emph{topology of fellow-traveling quasi-geodesics} and denote it $\mathcal{FQ}$. The idea is that a geodesic $\alpha$ is close to a geodesic $\beta$ if all quasi-geodesics tending to $\alpha$ closely fellow-travel quasi-geodesics tending to $\beta$ for a long time. 

One major difference from the direct limit topology on the Morse boundary is that rays with increasingly bad contracting functions can converge to a contracting ray. Recall again the space $\tilde{Y}$ from Example \ref{example:z fp zsqd} and the collection $\{\beta_i^j\}$ of rays from Remark \ref{rem: not 1st countable} . Consider the sequence $\{\beta^i_i\}$. It is not hard to see that this converges in the topology of the fellow-traveling quasi-geodesics to the geodesic ray $\alpha=aaaa\ldots$. The set $\{\beta_i^i\} \cup \{\alpha\}$ will be compact in $\mathcal{FQ}$. We note that the ray $\beta_i^i$ has a constant contracting function $\rho_i= i$, so this compact set has arbitrarily bad contracting geodesics.

The $\mathcal{FQ}$ topology keeps many of the desirable properties of the Morse boundary with the direct limit topology.  First they show  that it is a quasi-isometry invariant. Second, they also show that this boundary is second countable, and thus metrizable. Finally, they prove a weak version of North-South dynamics for the action of a group on its contracting boundary in \`a la Theorem \ref{thm: weak n-s} by Murray. Third, they show that if you restrict the $\mathcal{FQ}$ topology to rays which live in a single stratum and take the direct limit then this is homeomorphic to the Morse boundary. Finally they show that if $X$ is hyperbolic then this boundary is homeomorphic to the Gromov boundary. In summary:

\begin{theorem}[\cite{Cashen:2017aa}] Let $X$ be a proper geodesic space. The contracting boundary with the $\mathcal{FQ}$ topology, $\partial_c^\mathcal{FQ}$, is
\begin{enumerate}
\item a quasi-isometry invariant;
\item second countable and thus metrizable;
\item $\varinjlim_\rho \partial^\mathcal{FQ}_c X |_{\rho-\text{contracting}}$ is homeomorphic to the Morse boundary;
\item homeomorphic to the Gromov boundary if $X$ is hyperbolic;
\item has weak North-South dynamics  \`a la Theorem \ref{thm: weak n-s}.
\end{enumerate}
\end{theorem}

There are still many open questions about this topology: It is know that this topology is not in general homeomorphic to the subspace topology, but the example given is not a space with a geometric groups action. Is this topology different from the subspace topology in the presence of a geometric group action? We know the space is metrizable, but can we give a useful description of a metric? If so can we show that if $q \colon X \to Y$ is a quasi-isometry then $\partial q$ a quasi-symmetry?

\bibliographystyle{alpha}
\bibliography{../../../Library/library.bib}

\end{document}